\documentclass[reqno]{amsart}

\usepackage{bm}
\usepackage{amsmath}
\usepackage{amsfonts}
\usepackage{amssymb}
\usepackage{amsthm}
\usepackage{graphicx}
\usepackage{color}

\theoremstyle{plain}

\newtheorem{theorem}{Theorem}[section]
\newtheorem{corollary}[theorem]{Corollary}
\newtheorem{definition}[theorem]{Definition}
\newtheorem{example}[theorem]{Example}
\newtheorem{lemma}[theorem]{Lemma}
\newtheorem{notation}[theorem]{Notation}
\newtheorem{proposition}[theorem]{Proposition}
\newtheorem{remark}[theorem]{Remark}

\numberwithin{equation}{section}

\newcommand{\ordp}{\mathop{\mathrm{ord_p}}}
\newcommand{\ord}{\mathop{\mathrm{ord}}}

\begin{document}
\title[On the functional equation $Af^2+Bg^2=1$]{On the functional equation $Af^2+Bg^2=1$ on the field of complex $p$-adic numbers}\thanks{Work supported by FWF grant P16742-N04}
\author[E.~Mayerhofer]{Eberhard Mayerhofer}
\address{Vienna Institute of Finance, University of Vienna and Vienna University of Economics and Business Administration, Heiligenst\"adterstrasse 46-48, 1190 Vienna, Austria}
\email{eberhard.mayerhofer@vif.ac.at}
\begin{abstract}
For a fixed prime $p$, let $\mathbb C_p$ denote the complex $p$-adic numbers. For polynomials $A,B\in \mathbb C_p[x]$ we consider decompositions $A(x)f^2(x)+B(x)g^2(x)=1$ of entire functions $f,\,g$ on $\mathbb C_p$ and try to improve an impossibility result due to A. Boutabaa concerning transcendental $f,g$. We also provide an independent proof of a $p$-adic diophantic statement due to D. N. Clark, which is an important ingredient of Boutabaa's method.
\end{abstract}
\keywords{$p$-adic Analysis, $p$-adic functional equations, $p$-adic differential equations}
\subjclass[2000]{Primary: 39B52; Secondary: 12H25, 12J25, 30806 }
\maketitle
\section{Introduction}
For a prime $p$, let $\mathbb C_p$ denote the field of complex $p$-adic numbers, and let $\overline{\mathbb Q}\subset\mathbb C_p$ be the field of algebraic numbers over $\mathbb Q$. Furthermore, we denote by
$\mathcal A(\mathbb C_p)$ the ring of entire functions on $\mathbb C_p$. For a field $K$ we denote by
$K(x)$ the ring of rational functions over $K$, and $K[x]$ the ring of polynomials.

A. Boutabaa \cite{zB1} has shown:
\begin{theorem}\label{theoremB}
Let $A,\;B\in\overline{\mathbb Q}(x)$ be not identically zero. Suppose that $f,\,g\in\mathcal A(\mathbb C_p)$ have coefficients in $\overline{\mathbb Q}$ and satisfy
\begin{equation}\label{Af2}
A(x)f^2(x)+B(x)g^2(x)=1.
\end{equation}
Then $f,\,g\in\mathbb C_p[x]$.
\end{theorem}

The main aim of this paper is to provide a self-contained proof of the following more general result in the case
that $A, \,B$ are polynomials:
\begin{theorem}\label{theorem23}
Let $A,\;B\in\overline{\mathbb Q}[x]$ be not identically zero. Suppose that $f,\,g\in\mathcal A(\mathbb C_p)$
satisfy \eqref{Af2} and $g(0)\neq 0$. If $f^{(i)}(0)\in\overline{\mathbb Q},\;(0\leq i\leq\frac{\deg A+\deg B}{2}-1)$ and
 $g^{(i)}(0)\in\overline{\mathbb Q},\;(0\leq i<\frac{\deg A+\deg B}{2}-1)$, then $f,\,g\in\mathbb C_p[x]$.
\end{theorem}
When $\deg A$ and $\deg B$
do not have the same parity, we have the following stronger conclusion:
\begin{remark}\label{remark22}
Let $A$, $B$ be in $\mathbb C_p[x]$ such that $\deg A\not\equiv\deg B\mod 2$. Then equation \eqref{Af2} has no solutions in $(\mathcal A(\mathbb C_p))^2\setminus \mathbb C_p^2$.
\end{remark}
We note that the method of this paper allows for a similar assertion as that of Theorem \ref{theorem23} in the more abstract setting, when $\mathbb C_p$ is replaced by an algebraically closed and topologically complete ultrametric field $K$ of characteristic zero. This purely academic generalization however does not give more insight into the problem but only complicates notation and proofs. On the other hand, the case of $K$ with non-zero characteristic clearly makes no sense here.

\subsection{History}
The problem of decomposing complex meromorphic functions $f, g$ in the form $A(x)f^n+B(x)g^m=1$, where $A, B$ are certain meromorphic coefficients, has been thoroughly studied in the sixties and seventies of the last century. F.\ Gross \cite{zFG} shows that $f^n+g^n=1$ has no meromorphic solutions, when $n>3$ and no entire
solutions, if $n>2$. Moreover, for $n=2$ he characterizes all meromorphic solutions $f, g$.
Generalizations of these results are due to C.-C. Yang \cite{zCC} and N.\ Toda \cite{zNT}, and they are essentially applications of Nevanlinna Theory \cite{zIL}, in particular of the second Nevanlinna Theorem.
To see the match with our work, let us specialize C.-C. Yang's quite general result \cite{zCC} as follows:
\begin{theorem} \label{th yang}For non-constant complex entire functions $f, g$, the functional equation
\begin{equation}\label{Yang}
A(x)f^m(x)+B(x)g^n(x)=1,
\end{equation}
where $A,\,B\in \mathbb C[x]$ cannot hold, unless $n\leq 2,m\leq 3$.
\end{theorem}
Note that in the case $n=m=2$, $A=B=1$, a well known pair of entire solutions are the sine and cosine.

In the $p$-adic domain, A. Boutabaa considers decompositions of this kind in \cite{zB1}.
Building on an improved version of the $p$-adic second Nevanlinna Theorem \cite{zBE} it could be shown that compared with Theorem \ref{th yang}, the respective $p$-adic results are stronger:
\begin{theorem}
For non-constant $p$-adic entire functions $f, g\in\mathcal A(\mathbb C_p)$, the functional equation
\eqref{Yang} cannot hold with $A,B\in \mathbb C_p[x]$, unless $n=m=2$.
\end{theorem}
Note, that the pair of $p$-adic functions $\sin(x), \cos(x)$ is not entire, see also  Example \ref{example42}. This shortcoming motivates A. Boutabaa's ``impossibility result'' Theorem \ref{theoremB} when $n=m=2$. Since the $p$-adic Nevanlinna Theory is not applicable in this case, Boutabaa employs a Diophantic approximation result for algebraic numbers due to D.N.Clark  \cite{zClark}.
\subsection{Program of paper}
We start by recalling some elementary facts on entire functions and diophantic approximations of algebraic numbers in section \ref{prep sec}. We then state a new proof of Clark's result (Proposition \ref{lemma33}, see also Remark \ref{remark setoyanagi}). In the final section \ref{sec main} we present the proof of Theorem \ref{theorem23}.
\section{Preliminaries}\label{prep sec}

\begin{notation}\rm
\begin{itemize}
\item $\vert\,\,\,\vert_p$ denotes the $p$-adic absolute value on the field of complex $p$-adic numbers $\mathbb C_p$, and the respective additive valuation is given by $\ord_p(x):=-\log_p(\vert x\vert_p)$, where $\log_p(\cdot)$ is the logarithm with base $p$.
\item The ``closed'' disk with center $a\in \mathbb C_p$ and radius $r>0$ equals
\[d(a,r):=\{x\in \mathbb C_p:\;\vert  x-a\vert\leq r\}.
\]

\item We denote by $\mathbb Z_p$ the ring of $p$-adic integers, which is the topological completion of $\mathbb N$ with respect to the metric induced by the $p$-adic absolute value $\vert\,\,\,\vert_p$.
\item For $\xi\in \mathbb C_p$ we define the hypergeometric coefficient $(\xi)_k$ inductively as
\[
(\xi)_1:=\xi,\quad (\xi)_k:=(\xi)_{k-1}(\xi-k+1),\quad \textrm{  when  } k>1.
\]
\item For a natural number $n$, $\sigma(n)$ is the sum of the digits in its expansion as a $p$-adic integer and the factorial satisfies
\begin{equation}\label{eq factorial}
\ordp(n!)=\frac{n-\sigma_p(n)}{p-1},
\end{equation}
see \cite[p. 241]{zAR}.
\item For a set
$\mathcal D$, $\vert \mathcal D\vert$ denotes its cardinality. Furthermore, for a real number $x$, $[x]$ denotes its integral part. Accordinngly, $\langle x\rangle:=x-[x]$ is the fractional part of $x$.
\item Finally, we denote by $\mathcal A(\mathbb C_p)$ the algebra of entire functions on $\mathbb C_p$, that is the set of formal power series $f(x)=\sum_{i\geq 0}c_ix^i$ with coefficients $c_i$ in $\mathbb C_p$ for any $i\geq 0$ such that either $\{i\in\mathbb N: c_i\neq 0\}$ is finite (that is $f\in \mathbb C_p[x]$), or such that $\forall \lambda \in\mathbb R: \lim_{n\rightarrow\infty} \ord_p(c_n)-\lambda n=\infty$.
\item $\overline Q$ is the algebraic closure of $\mathbb Q\subset\mathbb C_p$.
\end{itemize}
\end{notation}
It is well known, that for any $r>0$, $\mathcal A(\mathbb C_p)$ can be endowed by a multiplicative
ultrametric norm $\vert\vert\;\vert\vert(r)$ defined as
\[
\vert\vert f\vert\vert(r):=\max_{i\geq 0} \vert c_i\vert_p r^i=\max_{x\in d(0,r)} \vert f(x)\vert_p.
\]
Also, as in classical complex analysis, we have ``Liouville's Theorem'': If for $f\in\mathcal A(\mathbb C_p)$, $\vert\vert f\vert\vert(r)$ is bounded as $r\rightarrow\infty$, then $f\in \mathbb C_p$. Indeed, this is an immediate consequence of the following useful result on entire functions \cite{zHY}:
\begin{lemma}\label{entirefct}
For any $f\in \mathcal A(\mathbb C_p)$ there exists $\alpha_0\in \mathbb C_p$ and $r\in\mathbb N$ such that $f$ can be decomposed as
\[
\alpha_0 x^r\prod_{f(\alpha)=0}\left(1-\frac{x}{\alpha}\right),\;\alpha_0\in \mathbb C_p,\;r\in\mathbb N.
\]
\end{lemma}
By means of this result we can prove the assertion of Remark \ref{remark22}:
\begin{proof}
Let  $a_s,\ b_t$ be the leading coefficients of $A, B$
respectively and let $d(0, R)$ be a disk containing all zeros of $A$ and
$B$. It is well known, that when $\vert x\vert>R$, we have
$\vert A(x)\vert=\vert a_s\vert \vert x\vert^s, \vert B(x)\vert=\vert b_t\vert \vert x\vert^t$. Consider now $\Gamma:=\{ x\mid\;r_1<\vert x \vert <r_2\}$ with $r_2>r_1>R$ such that $f,\; g$
have no zero inside $\Gamma$. Let $k$ be the number of zeros of $f$
and let $l$ be the number of zeros of $g$ in
$d(0,r_1)$. So the number of zeros of $f^2 $ (resp.\ $g^2$) is $2k$ (resp.\
$2l$). Due to Lemma \ref{entirefct}, $\vert f(x)\vert$ is of the
form $ \vert \lambda x^k\vert$ inside $\Gamma$ and $\vert g(x)\vert$ is of the
form $\vert\mu x^l\vert$ inside $\Gamma$. Therefore $\vert A(x)f^2(x)\vert$ is of the form $\vert a_s\lambda
x^{s+2k}\vert$, $\vert B(x)g^2(x)\vert$ is of the form
$\vert b_t\mu x^{t+2l}\vert$. Since $Af^2+Bg^2=1$, the two functions $\vert a_s\lambda\vert r^{s+2k},\; \vert b_t\mu\vert r^{t+2l}$
must be equal in $]r_1,r_2[$. This contradicts that $s,t$ have different parity.
\end{proof}
The $p$-adic version of Liouville's Theorem  (and improvements of it, e.g.\ the $p$-adic Thue-Siegel-Roth Theorem, \cite{zR}) is well known:
\begin{lemma}\label{liouville}
For any $\alpha\in d(0,1)\cap\overline{\mathbb Q}\setminus\mathbb N$ there exists a constant $C$ such that
for all $n\in\mathbb N$,
\[
\ordp(\alpha-n)\leq C+k\log_p(n),
\]
where $k$ is the degree of $\alpha$ over $\mathbb Q$.
\end{lemma}
We prepare Clark's statement with the following elementary observation:
\begin{proposition}\label{sumx}
For any $s\in\mathbb N$ we have
\[
\lim_{n\rightarrow\infty}\frac{1}{N}\sum_{0\leq j\leq N,\,\ordp(j)\leq s}\ordp(j)=\frac{1-p^{-s}}{p-1}-\frac{s}{p^{s+1}}.
\]
\end{proposition}
\begin{proof}
For $i\in\mathbb N$ we introduce the sets $L_i^N:=\{k\leq N:p^i\mid k,p^{i+1}\nmid k\}$. Then the power of this set is $\vert L_i^N\vert=\left[\frac{N}{p^i}\right]-\left[\frac{N}{p^{i+1}}\right]$. For any pair $i\neq j$, $L_i^N$ and $L_j^N$ are disjoint, and
$\bigcup_{i=0}^s L_i^N=\{j\vert\;j\leq N\}$. Also note that for $k\in L_i^N$, we have $\ordp(k)=i$.
Therefore, due to $\sum_{j\leq N}\ordp(j)=\sum_{i=0}^{s}\vert L_i^N\vert(i)$ we have
\begin{align*}\nonumber
\lim_{n\rightarrow\infty}\frac{1}{N}\sum_{0\leq j\leq N}\ordp(j)&=\lim_{n\rightarrow\infty}\frac{1}{N}\sum_{i=0}^{s}\left(\frac{N}{p^i}-
\frac{N}{p^{i+1}}\right)(i)=\\&=\sum_{i=0}^{s}(\frac{1}{p^i}-\frac{1}{p^{i+1}})(i)=\frac{1-p^{-s}}{p-1}-\frac{s}{p^{s+1}}.
\end{align*}
\end{proof}
The following Lemma is essentially due to D. N. Clark (\cite{zClark}). However taking into account a comment by M. Setoyanagi (\cite {zMS}) on Clark's conclusions we find it advisable to include an independent proof of this result. Compared with the original formulation in \cite{zClark}, also numbers not in $d(0,1)$ are involved, and the notion 'non-Liouville number' is avoided.
\begin{lemma}\label{lemma33}
If $\alpha\in \mathbb C_p$ is algebraic over $\mathbb{Q}$ or $\alpha\notin \mathbb Z_p$, then for sufficiently large $m$
\begin{equation}\label{clarkf}
\lim_{N\rightarrow\infty}\frac{\sum_{i=m}^N \ordp(\alpha-i)}{N}=w(\alpha)
\end{equation}
where $w(\alpha)$, the so-called 'weight' of $\alpha$, is given by
\[
w(\alpha)=\begin{cases}
\frac{1-p^{-[r(\alpha)]}}{p-1}-\langle r(\alpha)\rangle p^{-[(r(\alpha)]-1}&: \quad \alpha\in d(0,1)\setminus\mathbb Z_p\\
\frac{1}{p-1}&:\quad \alpha\in\mathbb Z_p\\
\ordp(\alpha)&:\quad\alpha\not\in d(0,1)
\end{cases}
\]
with $r(\alpha):=\sup_{i\geq 0}\ordp(\alpha-i)$.
\end{lemma}
\begin{proof}
We distinguish 4 cases:
\\{\bf Case 1: $\alpha\in \mathbb N$.} With $m=\alpha+1$ we have $\ordp(\prod_{j=m}^N(\alpha-j))=\frac{N-m+1-\sigma(N-m+1)}{p-1}$. Since $\sigma(j)=O(\log(j))$ when $j\rightarrow\infty$, we obtain
\[
\frac{1}{N} \sum_{j=m}^N \ordp(\alpha-j)\rightarrow \frac{1}{p-1},
\]
as $N\rightarrow \infty$.
\\{\bf Case 2: $\alpha\in\mathbb Z_p\setminus\mathbb N$.} We first establish lower bounds for the sum (\ref{clarkf}). For every natural number $\beta>N$ one has
\[
\left\vert\frac{\beta(\beta-1)\dots(\beta-N)}{(N+1)!} \right\vert_p\leq 1.
\]
Since $\mathbb N$ is dense in $\mathbb Z_p$ and due to the ultrametric rule ''the strongest one wins", for any natural number $N$ we can find some $\alpha_N\in\mathbb N$ such that $\alpha_N>N$ and $\vert \alpha_N-i\vert_p=\vert \alpha-i\vert_p$, for all $0\leq i\leq N$.
Therefore,
\[
\left\vert\frac{\alpha(\alpha-1)\dots(\alpha-N)}{(N+1)!} \right\vert_p\leq 1
\]
and due to formula \eqref{eq factorial} we infer
\[
\lim\inf_{N\rightarrow\infty}\frac{1}{N} \sum_{j=m}^N \ordp(\alpha-j)\geq\frac{1}{p-1}.
\]
Upper bounds for this sum can be achieved as follows. By induction it holds that
\begin{equation}\label{zeq2}
\frac{(-1)^{n+1}}{\alpha(\alpha-1)\dots(\alpha-n)}=\sum_{i+j=n} \frac{(-1)^{j+1}}{i!j!(\alpha-j)}.
\end{equation}
Multiplying (\ref{zeq2}) by $n!$ and evaluating with respect to the norm $\vert\;\vert_p$ we conclude from $\vert {n  \choose i}\vert_p=\vert \frac{n!}{i!(n-i)!}\vert_p\leq 1$ that
\begin{equation}\label{zeq3}
\left\vert\frac{n!}{(\alpha)_{n+1}}\right\vert_p\leq \max_{0\leq j\leq n}\frac{1}{\vert \alpha-j\vert_p}.
\end{equation}
Due to Lemma \ref{liouville} there exists $C'\in\mathbb{R}^+$ such that for any natural number j, $\vert \alpha-j\vert_p\geq \frac{C'}{j^{k}}$ for $k=[\mathbb Q(\alpha):\mathbb Q]$. After inserting this inequality into (\ref{zeq3}) we take the logarithm to the base $p$ of the $n$-th root of (\ref{zeq3}). Monotonicity of the root yields
and formula \eqref{eq factorial} yield
\begin{equation}\label{zeq4}
\frac{\ordp((\alpha)_{n+1})}{n}\leq \frac{1}{p-1}+\frac{k \log_{p} n-\sigma_p(n)}{n}-\frac{\log_{p} C'}{n}
\end{equation}
and taking limits we have
\[
\overline{\lim}_{N\rightarrow\infty}\frac{\sum_{i=0}^N\ordp(\alpha-i)}{N}\leq \frac{1}{p-1}.
\]
\\{\bf Case 3: $\alpha\in d(0,1)\setminus\mathbb Z_p$.} This is the most tricky part. First, we observe that $r(\alpha)=\sup_{i\geq 0} \ordp (\alpha-i)=\max_{i'\in\mathbb Z_p} \ordp(\alpha-i')$ exists in $\mathbb R$, since $\mathbb Z_p$ is compactly contained in $\mathbb C_p$. Moreover, the maximum is taken on by a natural number, that is, there exists $m\in\mathbb N$ such that $r(\alpha)=\ordp(\alpha-m)$. By setting $\beta:=\frac{\alpha-m}{p^[r(\alpha)]}$, we obtain
\begin{align*}\nonumber
r(\beta)&=\sup_{j\geq 0} \ordp(\beta-j)=\sup \ordp(\frac{\alpha-m}{p^{[r]}}-j)=\\\nonumber
&=\sup_{j\geq 0}\ordp(\alpha-m-jp^{[r(\alpha)]})-[r(\alpha)]=\\\nonumber=&\ordp(\alpha-m)-[r(\alpha)]=r(\alpha)-[r(\alpha)]=\langle r(\alpha)\rangle.
\end{align*}
For a fixed natural number $j>m$ we have
\begin{align*}\\\nonumber
\ordp(\alpha-j)&=\ordp(\beta+\frac{m-j}{p^{[r(\alpha)]}})+[r(\alpha)]=[r(\alpha)]+\\\nonumber &+\begin{cases}
\langle r(\alpha)\rangle &,\quad \ordp(\frac{m-j}{p^{[r(\alpha)]}})\geq \langle r(\alpha)\rangle\\ \ordp(\frac{m-j}{p^{[r(\alpha)]}}),&\quad \ordp(\frac{m-j}{p^{[r(\alpha)]}})<\langle r(\alpha)\rangle
\end{cases}.
\end{align*}
The second case in the right hand side of the last formula follows from the ultrametric rule " the strongest one wins". Concerning the first one, we note that $r(\beta)=\ordp(\beta)=\langle r(\alpha)\rangle$, hence by the non-archimedean triangle's inequality,
$\ordp(\beta+\frac{m-j}{p^{[r(\alpha)]}})\geq \langle r(\alpha)\rangle$. The reversed inequality holds due to the maximality of $r(\beta)$.
Rewriting the latter formula we receive
\begin{equation}\\\nonumber
\ordp(\alpha-j)=\begin{cases}
r(\alpha), &\quad \ordp(j-m) \geq [r(\alpha)]+1\\
\ordp(j-m),&\quad \ordp(j-m) \leq [r(\alpha)]
\end{cases}.
\end{equation}
It remains to verify \eqref{clarkf}. Since $r(\alpha-m)=r(\alpha)$ we may assume without loss of generality $m=0$. With this choice we obtain
\begin{align*}
\sum_{j=1}^N \ordp(\alpha-j)&=r(\alpha)\vert\{j\leq N\vert \ordp(j)\geq [r(\alpha)]+1\}\vert\\&+\sum_{1\leq j\leq N,\;\ordp(j)\leq [r(\alpha)]}\ordp(j)\end{align*}
Note that $\vert\{j\leq N\vert \ordp(j)\geq [r(\alpha)]+1\}\vert=\vert\{j\leq N\vert p^{[r(\alpha)]+1}\mid j \}\vert=\left[\frac{N}{p^{[r(\alpha)]+1}}\right]$. Furthermore, due to Lemma \ref{sumx},
\[
\frac{1}{N}\sum_{1\leq j\leq N,\;\ordp(j)\leq [r(\alpha)]}\ordp(j)\rightarrow \frac{1-p^{-[r(\alpha)]}}{p-1}-\frac{[r(\alpha)]}{p^{[r(\alpha)]+1}}
\]
when $N\rightarrow\infty$. Therefore, we conclude that
\[
\lim_{N\rightarrow\infty}\frac{1}{N}\sum_{j\leq N}\ordp(\alpha-j)=\frac{1-p^{-{r[\alpha]}}}{p-1}-\langle r(\alpha)\rangle p^{-[r(\alpha)]-1},
\]
which finishes the proof in Case 3.
\\{\bf Case 4: $\alpha\not\in d(0,1)$.} This is the trivial case, since due to the ultrametric rule "the strongest one wins" and $\ordp(i)\geq 0$ for all $i\geq 0$, one has $\ordp(\alpha-i)=\ordp(\alpha)$ for all $i\geq 0$.
\end{proof}
The following remark is due:
\begin{remark}\label{remark setoyanagi}
The proof of  ``Case 3'' seems to be more transparent than the one given in (\cite{zClark}). Also, in ``Case 2'' no geometric argument using 'Newton Polygons' is involved, as has been in (\cite{zClark}). Instead, we decompose the left side term of (\ref{zeq2}) into partial fractions.
\end{remark}
An immediate consequence of Lemma \ref{lemma33} is the following
\begin{corollary}\label{corclark}
For any $P\in\overline{\mathbb Q}[x]\setminus\{0\}$ there exists an integer $m\geq 0$ such that
\[
\lim_{N\rightarrow\infty} \frac{1}{N}\sum_{i=m}^N \ordp(P(i))<\infty.
\]
\end{corollary}
\begin{proof}
Decompose $P$ into a product of linear terms $x-\alpha_j$, where $\alpha_j \in\overline{\mathbb Q} \;(j=1,\dots,\deg P)$ times a constant. Due to the additivity of the valuation we may apply Lemma \ref{lemma33} to each of the factors.
\end{proof}
\section{Proof of the main Theorem}\label{sec main}
\begin{lemma}\label{lemma41}
Let $f,\,g\in\mathcal A(\mathbb C_p)$ and $A\, B\in \mathbb C_p[x]\setminus\{0\}$ such that \eqref{Af2} holds. Then $f$ and $g$ satisfy the system of differential equations
\begin{equation}\label{zeq5}
A'f+2Af'=hg,\;\;B'g+2Bg'=-hf
\end{equation}
for a certain polynomial $h$ over $\mathbb C_p$ with $\deg h\leq \frac{\deg A+\deg B}{2}-1$. Furthermore, if $\deg A>-\infty$ (that is $A\neq 0$) we have for some $\gamma>0$, $\vert\vert A'f+2Af'\vert\vert(r)\leq \gamma r^{\deg A-1}\vert\vert f\vert\vert(r)$ when $r$ sufficiently large.\end{lemma}
\begin{proof}
Differentiating the identity $Af^2+Bg^2=1$ yields $f(A'f+2Af')=-g(B'g+2Bg')$. Since $f$ and $g$ have no common zeros, there exists $h\in\mathcal A(\mathbb C_p)$ due to Lemma \ref{entirefct} such that $A'f+2Af'=hg,\; B'g+2Bg'=-hf$. Moreover, due to the identity $Af^2+Bg^2=1$ we have for $r$ sufficiently large, $\vert\vert Af^2\vert\vert(r)=\vert\vert Bg^2\vert\vert(r)$, i.e.
$\frac{\vert\vert g\vert\vert(r)}{\vert\vert f\vert\vert(r)}=\gamma r^ {\frac{\deg A-\deg B}{2}}$ for some $\gamma>0$. Thus
\[
\vert \vert h\vert\vert(r)=\frac{\vert\vert (Bg^2)'\vert\vert(r)}{\vert\vert g\vert\vert(r)\vert\vert f\vert\vert(r)}\leq \gamma r^{\frac{\deg A+\deg B}{2}-1}
\]
so we see that $h$ is a polynomial of degree $\leq \frac{\deg A+\deg B}{2}-1$.\end{proof}
\begin{example}\label{example42}\rm
The equation $f^2+g^2=1$ has no non-constant entire solutions: For, if there exists $f$, $g$ s.t. $f^2+g^2=(f+ig)(f-ig)=1$, where $i=\sqrt{-1}$, then, by taking norms on both sides yields
$\vert\vert f+ig\vert\vert(r)\vert\vert f-ig\vert\vert(r)=1$. Now, if $f+ig\in\mathcal A(\mathbb C_p)\setminus\mathbb C_p$, then for $r\rightarrow \infty$ we have $\vert\vert f+ig\vert\vert(r)\rightarrow\infty$ which implies $\vert\vert f-ig\vert\vert(r)\rightarrow 0$, and therefore $f=ig$ which yields the contradiction $f^2+g^2=0$.\end{example}
\begin{example}\label{example43}\rm
For $a\in \mathbb C_p$, the equation $f^2+(x-a)g^2=1$ has no entire solutions but $f=\pm 1$, $g=0$: Assume $g\neq 0$. Due to Lemma \ref{lemma41}, the mentioned polynomial $h$ is identically zero, therefore $f'=0$, but $(x-a)g^2$ cannot be a constant.
\end{example}
\begin{lemma}\label{lemma45}
Let $f,g \in\mathcal{A}(\mathbb C_p)$ satisfy $A'f+2Af'=hg, B'g+2Bg'=-hf$ for $A, B,h\in \mathbb C_p[x],h\not\equiv 0$, then $f$ satisfies a linear differential equation of second order:
\begin{equation}\label{zeq6}
(4ABh)f''+(6A'Bh+2AB'h-4ABh')f'+(A'B'h+2A''Bh-2BA'h'+h^3)f=0
\end{equation}
\end{lemma}
\begin{proof}
Due to Lemma \ref{lemma41}, (\ref{zeq5}) holds for $f$ and $g$ which implies $g=\frac{1}{h}\left(A'f+2 A f'\right)$. Differentiation yields $g'=\frac{-h'}{h^2}\left(A'f+2 A f'\right)+\frac{1}{h}\left(A''f+3A'f'+2Af''\right)$. Inserting $g$, $g'$ into (\ref{zeq5}) yields (\ref{zeq6}).\end{proof}
\begin{definition}\label{definition47}\rm
Consider a differential equation of the following form:
\[
Q^{(2)}(x)y''+Q^{(1)}(x)y'+Q^{(0)}(x)y=0,\; Q^{(i)}\in \mathbb C_p[x]\; (i=0,1,2)\;\;(E)
\]
with $Q^{(2)}$ not identically zero. We define the characteristic number of (E) as $N(E)=\max_{i=0,1,2}{\deg Q^{(i)}-i}$, and the characteristic polynomial of (E) as
\[
P_E(\xi):=\sum_{\deg Q^{(j)}-j=N(E)}q^{(j)}_{\deg Q^{(j)}}(\xi)_j
\]
where $q^{(j)}_{\deg Q^{(j)}}$ ($j=0,1,2$) are the leading coefficients of the polynomials $Q^{(j)}\;(j=0,1,2)$.\end{definition}
The following remark highlights the meaning of Definition \ref{definition47}:
\begin{remark}\label{remark48}\rm
Assume that $f(x)=\sum_{i\geq 0}c_i x^i\in\mathcal A(\mathbb C_p)\setminus \mathbb C_p[x]$ solves (E). Then
there exists a corresponding recurrence relation for the coefficients of $f$:
\begin{equation}\label{zequation9}
P_t(n)c_{n+t}+\dots+P_1(n)c_{n+1}+P_E(n)c_n=0
\end{equation}
with certain polynomials $P_s\in \mathbb C_p[x]$, where $1\leq s\leq t$, $P_t\neq 0$ and where $P_E$ is the characteristic polynomial of (E).\end{remark}
\begin{notation}\label{notation49}\rm
We write the polynomials $A, B$ and $h$ (due to Lemma \ref{lemma41}) in the following way:
\[
A(x)=\sum_{i=0}^{\eta} a_ix^i,\;B(x)=\sum_{j=0}^{\chi}b_jx^j,\;h(x)=\sum_{k=0}^{\mu}h_kx^k,\;a_{\eta},b_{\chi}\neq0,
\]
i.e., $A, B$ are polynomials of degree $\eta$ resp. $\chi$. Note that $h$ might be identically zero.
\end{notation}
We now determine the characteristic polynomial of equation (\ref{zeq6}):
\begin{lemma}\label{Case2}
Let $A, B \in \overline{\mathbb Q}[x]$, $f,g\in\mathcal A(\mathbb C_p)$ such that \eqref{Af2} holds. Suppose the polynomial $h$ of Lemma (\ref{lemma41}) is not identically zero\footnote{this in particular implies $\deg A+\deg B\geq 2$} and that $g(0)\neq 0$,
$f^{(i)}(0)\in \overline{\mathbb Q}$ for $0\leq i\leq \frac{\deg A+\deg B}{2}-1$, and
$g^{(i)}(0)\in \overline{\mathbb Q}$ for $0\leq i< \frac{\deg A+\deg B}{2}-1$. Then the characteristic polynomial $P_E$ of the corresponding differential equation (\ref{zeq6}) is an element of $\overline{\mathbb Q}[x]$.
\end{lemma}
\begin{proof}
We make use of the notation introduced in Definition \ref{definition47} for the polynomial coefficients of (\ref{zeq6}):
$Q^{(2)}=4ABh$, $Q^{(1)}=6A'Bh+2AB'h-4ABh'$, $Q^{(0)}=A'B'h+2A''Bh-2BA'h'+h^3$. Recall that $a_{\eta},b_{\chi}$, $h_{\mu}$ denote the coefficients of the leading powers of the resp. polynomials. Note that
\[
N(E):=\max_{i=0,1,2} \deg Q{(i)}-i=\eta+\chi+\mu-2
\]
since $\deg Q^{(1)}\leq \eta+\chi+\mu-1$ and $\deg Q^{(0)}\leq \eta+\chi+\mu-2$. First, let us calculate
the coefficients $q_j$ of the term $x^{N(E)+j}$ in $Q^{(j)}\;(j=0,1,2)$ depending on the degree of $h$ ($q_j=0$ might vanish for $j=1$ or $j=0$!)
\begin{enumerate}
\item $\mu=\deg h=\frac{\deg A+\deg B}{2}-1=\frac{\eta+\chi}{2}-1$: Then $q_2=4a_{\eta}b_{\chi}h_{\mu}$, $q_1=a_{\eta}b_{\chi}h_{\mu}(6\eta+2\chi-4\mu)$ and finally $q_0=a_{\eta}b_{\chi}h_{\mu}(\eta\chi+2\eta(\eta-1)-2\eta\mu+\frac{h_{\mu}^2}{a_{\eta}b_{\chi}})$.
\item $\mu=\deg h<\frac{\deg A+\deg B}{2}-1=\frac{\eta+\chi}{2}-1$: Then, due to the calculations in the case above $q_2=4a_{\eta}b_{\chi}h_{\mu},\;q_1=a_{\eta}b_{\chi}h_{\mu}(6\eta+2\chi-4\mu)$ and $q_0=a_{\eta}b_{\chi}h_{\mu}(\eta\chi+2\eta(\eta-1)-2\eta\mu)$.
\end{enumerate}

We derive the characteristic polynomial $P_E(\xi)$:
\begin{enumerate}
\item $\mu=\deg h=\frac{\deg A+\deg B}{2}-1=\frac{\eta+\chi}{2}-1$: $P_E(\xi)=a_{\eta}b_{\chi}h_{\mu}[4 \xi(\xi-1)+4(\eta+1)\xi+(\eta^2+\frac{h_{\mu}^2}{a_{\eta}b_{\chi}})$.

\item $\mu=\deg h<\frac{\deg A+\deg B}{2}-1=\frac{\eta+\chi}{2}-1$: $P_E(\xi)=a_{\eta}b_{\chi}h_{\mu}[4\xi(\xi-1)+(6\eta+2\chi-4\mu)\xi+(\eta\chi+2\eta(\eta-1)-2\eta\mu)]$.
\end{enumerate}
Due to our assumptions it suffices to show that the leading coefficient $h_{\mu}$ of $h$ is algebraic over $\mathbb Q$. Indeed by differentiating $h=\frac{A'f+2Af'}{g}$ $\mu$-times and using Lemma \ref{lemma41} we infer
\[
\mu! h_{\mu}=(h(x))^{(\mu)}=\left(\frac{(A(x)f^2(x))'}{fg}\right)^{(\mu)}.
\]
The right hand side involves derivatives up to order $\mu+1$ of $f$
and of up to order $\mu$ of $g$. Due to our assumptions the coefficients of $f, g$ up to order $\mu+1$ resp. $\mu$ are algebraic over $\mathbb Q$, and since the coefficients of $A, \;B$ have the same property, we are done.
\end{proof}
We shall also employ the 'entire version' of \cite[Proposition]{zB2}):
\begin{lemma}\label{lemma411}
Consider the linear differential equation
\[
C(x)y'(x)+D(x)y(x)=0, \;\;C(x),D(x)\in K[x]
\]
with $C$ not identically zero. Let $y(x)\in\mathcal A(\mathbb C_p)$. Then $y(x)\in \mathbb C_p[x]$.
\end{lemma}
{\bf Proof of Theorem \ref{theorem23}}
Let $f, g\in \mathcal A(\mathbb C_p)\setminus \mathbb C_p[x]$ such that \eqref{Af2} holds. By Lemma \ref{lemma41}, $f,g$ satisfy equation (\ref{zeq5})
with a certain polynomial $h\in \mathbb C_p[x]$. We have the two cases:
\\\medskip {\bf Case 1: $h$ identically zero:}
Due to Lemma \ref{lemma411}, the solutions of system (\ref{zeq5}) are in $\mathbb C_p[x]$, which is impossible.
\\\medskip {\bf Case 2: $h$ not identically zero:}
Due to Lemma \ref{lemma45}, $f$ satisfies a linear differential equation of the form (\ref{zeq6}). Moreover, due to Lemma \ref{Case2}, the characteristic polynomial $P_E$ of (\ref{zeq6}) lies in $\overline{\mathbb Q}[x]$.
Following the notation of Remark \ref{remark48} we consider a recurrence relation for the coefficients $c_n$ of $f(x)=\sum_{n\geq 0}c_nx^n$, $n$ sufficiently large
\begin{equation}\label{eqp2}
P_t(n)c_{n+t}+\dots+P_1(n)c_{n+1}+P_E(n)c_n=0,\; t>0
\end{equation}
We consider the non-trivial case, where $P_t$ is not identically zero. Clearly
$P_E$ is not identically zero, because $h$ is not. Since we can multiply
equation (\ref{eqp2}) with an appropriate constant $\gamma$, we may assume without loss of generality that
$\forall i\in\{1,\dots,t\}\forall n\in\mathbb N:\ordp(P_i(n))\geq 0$ where we have
set $P_0:=\gamma P_E$. We have $P_{E}\in\gamma\overline{\mathbb Q}[x]$ due to Lemma \ref{Case2}. Moreover, by induction
it follows for any $k\in\mathbb N$ there exist certain polynomials
$P_{i,k}, i\in\{1,\dots,t\}$ with $\ordp(P_{i,k}(n))\geq 0\;\forall n\geq 1$ such that the following recurrence relation is satisfied for any $n$:
\begin{equation}\label{eqp3}
P_{t,k}(n)c_{n+t+k}+\dots+P_{1,k}(n)c_{n+k+1}+P_0(n)P_0(n+1)\dots P_0(n+k)c_n=0
\end{equation}
Evaluation with respect to $\ordp$ yields
\begin{multline}\label{eqp4}
\ordp(P_0(n)P_0(n+1)\dots P_0(n+k))+\ordp(c_n)\geq\\\geq\min\{\ordp(c_{n+t+k}),\dots,\ordp(c_{n+k+1})\}
\end{multline}
Moreover, due to Corollary \ref{corclark}, $\exists L\in \mathbb R:\lim_{k\rightarrow\infty}\frac{\ordp(P_0(n)P_0(n+1)\dots P_0(n+k))}{n+k}= L$. Consider now the growth of $c_n$, the coefficients of $f$: the transcendence of $f$ implies that $ \forall \lambda \in\mathbb R: \lim_{n\rightarrow \infty} \ordp(c_n)-\lambda n=\infty$. In other words, for any $\lambda\in\mathbb R$ and for any real constant $c(\lambda)$ we have $\ordp(c_n)\geq c(\lambda)+\lambda n $ for sufficiently large $n$.
We may choose $\lambda >L$ and some $c(\lambda)$. Applying this to (\ref{eqp4})
divided by $n+k$ yields for sufficiently large $n$:
\begin{equation}\label{eqp6}
\frac{\ordp(\delta^kP_0(n)P_0(n+1)\dots P_0(n+k))}{n+k}+\frac{\ordp(c_n)}{n+k}\geq \frac{c(\lambda)+\lambda (n+k+1)}{n+k}.
\end{equation}
Taking the limites on both sides ($k \rightarrow\infty)$ we derive the contradiction $L\geq \lambda$.\\
\\{\bf Acknowledgement:}
For reading this manuscript carefully I am indebted to Professor Alain Escassut (Clermont-Ferrand, France). Special thanks to the referee for suggesting a number of simplifications which lead to the final form of this paper.

\end{document}